\documentclass[reqno, 12pt]{amsart}
\pdfoutput=1
\makeatletter
\let\origsection=\section \def\section{\@ifstar{\origsection*}{\mysection}}
\def\mysection{\@startsection{section}{1}\z@{.7\linespacing\@plus\linespacing}{.5\linespacing}{\normalfont\scshape\centering\S}}
\makeatother

\usepackage{amsmath,amssymb,amsthm}
\usepackage{mathrsfs}
\usepackage{mathabx}\changenotsign
\usepackage{dsfont}

\usepackage[normalem]{ulem}

\usepackage[dvipsnames]{xcolor}
\usepackage{hyperref}
\hypersetup{
    colorlinks,
    linkcolor={red!60!black},
    citecolor={green!60!black},
    urlcolor={blue!60!black}
}

\usepackage{graphicx}

\usepackage[open,openlevel=2,atend]{bookmark}

\usepackage[abbrev,msc-links]{amsrefs}
\usepackage{doi}

\renewcommand{\PrintDOI}[1]{\doi{#1}}

\usepackage[T1]{fontenc}
\usepackage{lmodern}
\usepackage[babel]{microtype}
\usepackage[english]{babel}

\usepackage{geometry}
\numberwithin{equation}{section}
\numberwithin{figure}{section}

\usepackage{enumitem}

\usepackage{pgf,tikz}

\makeatletter
\def\greek#1{\expandafter\@greek\csname c@#1\endcsname}
\def\Greek#1{\expandafter\@Greek\csname c@#1\endcsname}
\def\@greek#1{\ifcase#1
	\or $\alpha$%
	\or $\beta$%
	\or $\gamma$%
	\or $\delta$%
	\or $\epsilon$%
	\or $\zeta$%
	\or $\eta$%
	\or $\theta$%
	\or $\iota$%
	\or $\kappa$%
	\or $\lambda$%
	\or $\mu$%
	\or $\nu$%
	\or $\xi$%
	\or $o$%
	\or $\pi$%
	\or $\rho$%
	\or $\sigma$%
	\or $\tau$%
	\or $\upsilon$%
	\or $\phi$%
	\or $\chi$%
	\or $\psi$%
	\or $\omega$%
\fi}
\def\@Greek#1{\ifcase#1
	\or $\mathrm{A}$%
	\or $\mathrm{B}$%
	\or $\Gamma$%
	\or $\Delta$%
	\or $\mathrm{E}$%
	\or $\mathrm{Z}$%
	\or $\mathrm{H}$%
	\or $\Theta$%
	\or $\mathrm{I}$%
	\or $\mathrm{K}$%
	\or $\Lambda$%
	\or $\mathrm{M}$%
	\or $\mathrm{N}$%
	\or $\Xi$%
	\or $\mathrm{O}$%
	\or $\Pi$%
	\or $\mathrm{P}$%
	\or $\Sigma$%
	\or $\mathrm{T}$%
	\or $\mathrm{Y}$%
	\or $\Phi$%
	\or $\mathrm{X}$%
	\or $\Psi$%
	\or $\Omega$%
\fi}

\AddEnumerateCounter{\greek}{\@greek}{24}
\AddEnumerateCounter{\Greek}{\@Greek}{12}

\makeatother

\let\polishlcross=\l
\def\l{\ifmmode\ell\else\polishlcross\fi}

\def\paragraph#1{%
  \noindent\textbf{#1.}\enspace}

\makeatletter
\def\moverlay{\mathpalette\mov@rlay}
\def\mov@rlay#1#2{\leavevmode\vtop{   \baselineskip\z@skip \lineskiplimit-\maxdimen
   \ialign{\hfil$\m@th#1##$\hfil\cr#2\crcr}}}
\newcommand{\charfusion}[3][\mathord]{
    #1{\ifx#1\mathop\vphantom{#2}\fi
        \mathpalette\mov@rlay{#2\cr#3}
      }
    \ifx#1\mathop\expandafter\displaylimits\fi}
\makeatother

\DeclareFontFamily{U}  {MnSymbolC}{}
\DeclareSymbolFont{MnSyC}         {U}  {MnSymbolC}{m}{n}
\DeclareFontShape{U}{MnSymbolC}{m}{n}{
    <-6>  MnSymbolC5
   <6-7>  MnSymbolC6
   <7-8>  MnSymbolC7
   <8-9>  MnSymbolC8
   <9-10> MnSymbolC9
  <10-12> MnSymbolC10
  <12->   MnSymbolC12}{}
\DeclareMathSymbol{\powerset}{\mathord}{MnSyC}{180}

\usepackage{tikz}
\usetikzlibrary{calc,decorations.pathmorphing}
\pgfdeclarelayer{background}
\pgfdeclarelayer{foreground}
\pgfdeclarelayer{front}
\pgfsetlayers{background,main,foreground,front}

\let\epsilon=\varepsilon

\let\rho=\varrho
\let\theta=\vartheta
\let\kappa=\varkappa

\def\EE{{\mathds E}}
\let\E=\EE

\def\PP{{\mathds P}}

\theoremstyle{plain}
\newtheorem{thm}{Theorem}[section]
\newtheorem{theorem}[thm]{Theorem}

\newtheorem{lemma}[thm]{Lemma}

\theoremstyle{definition}

\newtheorem{exmp}[thm]{Example}
\newtheorem{conj}[thm]{Conjecture}

\usepackage{accents}

\let\phi=\varphi

\begin{document}

\title[Multiple twins in permutations]{Multiple twins in permutations}

\author{Andrzej Dudek}
\address{Department of Mathematics, Western Michigan University, Kalamazoo, MI, USA}
\email{\tt andrzej.dudek@wmich.edu}
\thanks{The first author was supported in part by Simons Foundation Grant \#522400.}

\author{Jaros\l aw Grytczuk}
\address{Faculty of Mathematics and Information Science, Warsaw University of Technology, Warsaw, Poland}
\email{j.grytczuk@mini.pw.edu.pl}
\thanks{The second author was supported in part by Narodowe Centrum Nauki, grant 2017/26/D/ST6/00264, as well as, the European Regional Development Fund under the grant No. POIR.01.01.01-00-0124/17-00, on the basis of an agreement between FinAi S.A. and the National Center for Research and Development based in Warsaw.}

\author{Andrzej Ruci\'nski}
\address{Department of Discrete Mathematics, Adam Mickiewicz University, Pozna\'n, Poland}
\email{\tt rucinski@amu.edu.pl}
\thanks{The third author was supported in part by Narodowe Centrum Nauki, grant 2018/29/B/ST1/00426}

\begin{abstract}
By an $r$-tuplet in a permutation we mean a family of $r$ pairwise disjoint subsequences with the same relative order. The length of an $r$-tuplet is defined as the length of any single subsequence in the family. Let $t^{(r)}(n)$ denote the largest $k$ such that every permutation of length $n$ contains an $r$-tuplet of length $k$. We prove that $t^{(r)}(n)=O\left(n^{\frac r{2r-1}}\right)$ and $t^{(r)}(n)=\Omega\left( n^{\frac{R}{2R-1}} \right)$, where $R=\binom{2r-1}r$. We conjecture that the upper bound brings the correct order of magnitude of  $t^{(r)}(n)$ and support this conjecture by proving that it holds for almost all permutations. Our work generalizes previous studies of the case $r=2$.

\end{abstract}

\maketitle


\setcounter{footnote}{1}

\section{Introduction}

By a \emph{permutation} we mean  any finite sequence of  \emph{distinct} integers. We say that two permutations $(x_1,\dots,x_k)$ and $(y_1,\dots,y_k)$ are \emph{similar} if their entries preserve the same relative order, that is, $x_i<x_j$ if and only if $y_i<y_j$ for all pairs $\{i,j\}$ with $1\leqslant i<j\leqslant k$. This is, clearly, an equivalence relation. Note that given a permutation $(x_1,\dots,x_k)$ and a $k$-element set $\{y_1,\dots,y_k\}$ of positive integers, there is only one permutation of this set similar to $(x_1,\dots,x_k)$.

Let $[n]=\{1,2,\dots,n\}$ and $\pi$ be a permutation of $[n]$, called also an \emph{$n$-permutation}, and let $r\ge2$ be an integer. A family of $r$ pairwise similar and disjoint sub-permutations of $\pi$ is called an \emph{$r$-tuplet} and the \emph{length} of an $r$-tuplet is defined as the number of elements in just one of the sub-permutations. It is common to call $2$-tuplets \emph{twins} and $3$-tuplets \emph{triplets}.

Let $t^{(r)}(\pi)$ denote the largest integer $k$ such that $\pi$ contains an $r$-tuplet of length $k$. Further, let $t^{(r)}(n)$ denote the minimum of $t^{(r)}(\pi)$ over all permutations $\pi$ of $[n]$. In other words, $t^{(r)}(n)$ is the largest integer $k$ such that every $n$-permutation contains an $r$-tuplet of length $k$. Our aim is to estimate this function.

By the classical result of Erd\H {o}s and Szekeres \cite{ErdosSzekeres} concerning monotone subsequences, we get $t^{(r)}(n)= \Omega(\sqrt{n})$ (all implicit constants throughout the paper are allowed to depend on $r$). Indeed, in a monotone permutation any $r$ disjoint subsequences of the same length yield an $r$-tuplet. For $r=2$, using a probabilistic argument, Gawron \cite{Gawron} proved that $t^{(2)}(n)=O(n^{2/3})$. He  conjectured that this bound is tight, namely, that we also have $t^{(2)}(n)=\Omega(n^{2/3})$. A currently best result towards this conjecture was obtained by Bukh and Rudenko~\cite{BukhR} (see also~\cite{DGR}) who showed that $t^{(2)}(n)=\Omega(n^{3/5})$.

In the first part of this paper, we generalize both these bounds to arbitrary $r\ge2$, adopting the ideas from  \cite{BukhR, DGR, Gawron}.

\begin{thm}\label{thm:_ulbounds} For every $r\ge2$, with $R=\binom{2r-1}r$, we have
\[
t^{(r)}(n) = \Omega\left( n^{\frac{R}{2R-1}} \right)  \quad \text{ and } \quad t^{(r)}(n) =  O\left(n^{\frac r{2r-1}}\right).
\]
\end{thm}
\noindent

Let $\Pi_n$ be a (uniformly) random permutation of $[n]$ and let $t^{(r)}(\Pi_n)$ be the corresponding random variable equal to the maximum length of an $r$-tuplet in $\Pi_n$. We say that a property of a random object holds \emph{asymptotically almost surely} (a.a.s.~for short) if it holds with probability tending to one as the size of the object grows to infinity. In \cite{DGR} and \cite{BukhR} it was shown that a.a.s., $t^{(2)}(\Pi_n)=\Theta(n^{2/3})$. Here we generalize this result.

\begin{thm}\label{rg1}
For every $r\ge2$, a.a.s.,
\[
t^{(r)}(\Pi_n) = \Theta\left(n^{\frac{r}{2r-1}}\right).
\]
\end{thm}

\noindent
In view of Theorem \ref{rg1}, for almost all permutations $\pi$ the parameter $t^{(r)}(\pi)$ reaches the upper bound from Theorem~\ref{thm:_ulbounds}.

In the following two sections we give proofs of Theorems \ref{thm:_ulbounds} and \ref{rg1}, while the last section contains some related open problems.

\section{Proof of Theorem \ref{thm:_ulbounds}}

The upper bound in Theorem~\ref{thm:_ulbounds} follows from the upper bound in Theorem~\ref{rg1}. Therefore, here we focus exclusively on
the proof of the lower bound.

\subsection{Preparations}\label{subsec:prep}

Beame and Huynh-Ngoc~\cite{BH-N}  proved that amongst any three permutations of $[m]$ there are two with the same sub-permutation of length at least $m^{1/3}$. This simple fact was used in \cite{BukhR} to prove the lower bound in Theorem \ref{thm:_ulbounds} for $r=2$. Here, we will need the following extension  proved by Bukh and Zhou \cite{BZ} in the context of estimating the length of $r$-tuplets in words over finite alphabets.
Recall that $R=\binom{2r-1}r$.

\begin{lemma}[Theorem 24 in \cite{BZ}]\label{3to2r-1}
For every $r\ge 2$, among any $2r-1$ permutations of $[m]$, there are $r$ permutations with the same sub-permutation of length at least $m^{1/R}$.
\end{lemma}

The proof of Lemma \ref{3to2r-1} is quite elementary and uses the idea from the classical proof of the Erd\H os-Szekeres Theorem. It boils down to assigning to each $i\in[m]$ a vector of length $R$, where the entries represent, for all $r$-tuples of permutations, the lengths of the longest common sub-permutations  that begin at $i$. Then it suffices to show that this mapping is an injection.

\bigskip

The main technique behind the proof of the lower bound on $t^{(r)}(n)$ is that of \emph{concatenation} of $r$-tuplets. Here we reveal sufficient conditions under which two $r$-tuplets can be merged into one. Let $(x^{(j)}_{1},\dots,x^{(j)}_{s})$, $j=1,\dots,r$, and $(y^{(j)}_{1},\dots,y^{(j)}_{t})$, $j=1,\dots,r$, be two disjoint $r$-tuplets in a permutation $\pi$. Their full concatenation, resulting in an $r$-tuplet of length $s+t$, is possible if
\begin{itemize}
 \item the rightmost element of the first $r$-tuplet is to the left of the leftmost element of the second $r$-tuplet, and
  \item for each $1\le i\le s$ and $1\le i'\le t$,
  \[ \text{either} \qquad  \max_{1\le j\le r}y^{(j)}_{i'}<\min_{1\le j\le r}x^{(j)}_{i} \qquad \text{ or } \qquad \min_{1\le j\le r}y^{(j)}_{i'}>\max_{1\le j\le r}x^{(j)}_{i}.\]
\end{itemize}
\begin{exmp} Let $r=2$, $n=30$, and
   \begin{align*}\pi=(&\underline{\colorbox{cyan}{26}},16,28,29,\underline{\colorbox{cyan}{10}},5,\overline{\colorbox{pink}{24}},\underline{\colorbox{cyan}{27}},
   1,22,11,\overline{\colorbox{pink}{8}},2,23,15,
   \\&19,\overline{\colorbox{pink}{25}},21,20,\underline{\underline{\colorbox{RoyalBlue}{13}}},9,30,\underline{\underline{\colorbox{RoyalBlue}{17}}},
  \overline{\overline{\colorbox{red}{12}}},\overline{\overline{\colorbox{red}{18}}},7,\underline{\underline{\colorbox{RoyalBlue}{3}}},14,\overline{\overline{\colorbox{red}{4}}},6)
   \end{align*}
   be a permutation of $[30]$. 
Here we marked by blue (single underline) and pink (single overline) the first pair of twins of length 3 (similar to $(2,1,3)$), and by indigo (double underline) and red (double overline) the second one (similar to $(2,3,1)$). Both above conditions hold. Indeed,  elements of the first twins are to the left of 13, the leftmost element of the second twins. Moreover, all nine required inequalities hold, e.g., for $i=2$, $i'=1$, we have $\min\{13,12\}=12>\max\{10,8\}=10$. Thus,  their concatenation  forms a (bluish-reddish) pair of twins of length six, namely: $(26, 10,27,13,17,3)$ and $(24,8,25,12,18,4)$, both similar to $(5,2,6,3,4,1)$.

If we only knew the first twins and were after the second ones, we could, obviously, help ourselves by searching only to the right of 25, and by eliminating from our search all elements whose values are ``squeezed'' between the given twins, that is, lie between $\min_{1\le j\le r}x^{(j)}_{i}$ and $\max_{1\le j\le r}x^{(j)}_{i}$ for any $i$. (In our example this step  excludes only 9 ($i=2$), as 25 ($i=1$) and 26 ($i=3$) are already excluded.) But most importantly, the second pair of twins should be ``narrow'' in a sense that would guarantee the second bulleted condition above. (In our example, we chose consecutive values: $13,12$ and $17,18$ and $3,4$.)
\end{exmp}

To facilitate the idea mentioned in the above example, we introduce the notion of narrow $r$-tuplets as follows. The \emph{width} of a set of integers $A$ is defined as $\max A-\min A$. For a positive integer $w$ we call an $r$-tuplet $(x^{(1)}_1,\dots,x^{(1)}_k),\dots,(x^{(r)}_1,\dots,x^{(r)}_k)$ \emph{$w$-narrow} if for all $i=1,\dots,k$, the sets $\{x^{(j)}_i: 1\le j\le r\}$ have width at most $w$.

Throughout the paper we will sometimes use notation $[a,b]=\{a,\dots,b\}$, where $a<b$ are integers.

\subsection{Lower bound}

   The proof of the lower bound is similar to that for $r=2$ in \cite[Proof of Theorem 1]{BukhR}. The main idea  is to utilize Lemma \ref{3to2r-1} in a clever way.

   Set
   $$M=n^{(R-1)/(2R-1)}\quad\mbox{and}\quad N=n/M=n^{R/(2R-1)}.$$
   Our goal is, for a given permutation $\pi$ of $[n]$, to gradually pick from it $O(M)$-narrow $r$-tuplets of length $\Omega(N/M)$, and remove them together with all ``in-between'' elements (plus some more), allowing concatenation of obtained pieces into one  $r$-tuplet of length $\Omega(N)$. To achieve this goal, we need to iterate this procedure $\Omega(M)$ times, meaning that the number of elements discarded each time  should not exceed $O(N)$. Now come the details.

     Let $A\subset\binom{[n]}{(2r-1)N}$, $A=\{a_1<a_2\cdots<a_{(2r-1)N}\}$, and let $\pi$ be a permutation of $A$. For any $a\in A$ we denote by $\pi^{-1}(a)$ the position of $a$ in $\pi$. E.g., if $A=\{1,3,4,6\}$ and $\pi=(3,1,6,4)$, we have $\pi^{-1}(1)=2$, $\pi^{-1}(3)=1$, $\pi^{-1}(4)=4$ and $\pi^{-1}(6)=3$.

     We split the set $A$ into $N$ blocks of $2r-1$ consecutive elements:
     $$A=A_0\cup\cdots\cup A_{N-1},\quad\mbox{where}\quad A_i=\{a_{i(2r-1)+1},\dots,a_{(i+1)(2r-1)}\},\quad i=0,\dots,N-1.$$
   A crucial observation is that if, for some $2\le k\le N$ and $1\le j_1<j_2\le 2r-1$, there exists a sequence of \emph{distinct} indices $(i_1,\dots,i_k)\in\{0,\dots,N-1\}^k$ such that
   $$\pi^{-1}(a_{i_1(2r-1)+j_1})<\cdots<\pi^{-1}(a_{i_k(2r-1)+j_1})\quad\mbox{and}\quad\pi^{-1}(a_{i_1(2r-1)+j_2})<\cdots<\pi^{-1}(a_{i_k(2r-1)+j_2}),$$
   then sub-permutations  $(a_{i_1(2r-1)+j_1},\dots,a_{i_k(2r-1)+j_1})$ and $(a_{i_1(2r-1)+j_2},\dots,a_{i_k(2r-1)+j_2})$ of~$\pi$ form twins of length~$k$. By the same token, if all the above is true for $r$ indices $1\le j_1<j_2<\cdots < j_r\le 2r-1$, then we obtain an $r$-tuplet of length $k$ in~$\pi$.

   \begin{exmp} Let $A=[30]$, $r=3$, $N=6$, and, same as before,
   \begin{align*}
   \pi=(&26,16,28,29,10,5,24,27,1,22,11,8,2,23,15,\\&19,25,21,20,13,9,30,17,12,18,7,3,14,4,6).
   \end{align*}
   Here $A_0=[1,5]$, $A_1=[6,10]$, $A_2=[11,15]$, $A_3=[16,20]$, $A_4=[21,25]$ and $A_5=[26,30]$.
   By inspection, one can see that for $k=3$, $i_1=5$, $i_2=3$, $i_3=0$, and $j_1=1$, $j_2=3$, $j_3=4$, the above conditions are satisfied.
   Indeed, the first element of $A_5$ (26 in $\pi$) is to the left of the first element of $A_3$ (16) which, in turn, is to the left of the first element of $A_0$ (1), and similar order holds for the third elements of $A_5$, $A_{3}$, $A_0$ (28, 18, 3), as well as, for the fourth elements (29, 19, 4).
   This reveals a triplet similar to $(3,2,1)$ and indicated by colors blue (underline), pink (overline), and red (asterisk)  below:
   \begin{align*}\pi=(&\underline{\colorbox{cyan}{26}},\underline{\colorbox{cyan}{16}},\overline{\colorbox{pink}{28}},\colorbox{red}{29*},10,5,24,27,
   \underline{\colorbox{cyan}{1}},22,11,8,2,23,15,\\
   &\colorbox{red}{19*},25,21,20,13,9,30,17,12,\overline{\colorbox{pink}{18}},7,
   \overline{\colorbox{pink}{3}},14,\colorbox{red}{4*},6).
   \end{align*}

   An equivalent but somewhat easier way to see what is going on here is to define five disjoint sub-permutations $\pi^{(j)}$ of length 6, by including in $\pi^{(j)}$ all elements of $\pi$ of the form $a_{i(2r-1)+j}$ (in the order they appear in $\pi$): $\pi^{(1)}=(26,16,1,11,21,6)$, $\pi^{(2)}=(27,22,2,17,12,7)$, $\dots$, etc. Next replace them by similar permutations of $\{0,1,2,3,4,5\}$: $(5,3,0,2,4,1)$, $(5,4,0,3,2,1)$, $\dots$, etc. Now, what we are after, are long common sub-permutations in at least three of these permutations. We found $(5,3,0)$ as common in the 1st, 3rd, and 4th permutation, but one could have also picked $(5,3,1)$ in the 1st, 2nd, and 4th (and there is at least one more alternative). By inspection, one may realize, however, that there are no longer common sub-permutations in at least three of these five permutations.
   \end{exmp}

In the above example, we dealt with consecutive integers in $A$, so the width of each subset $A_i$ was the same (and equal to 4), but in general this may not be the case. To make sure that the $r$-tuplets generated by the above approach are $O(M)$-narrow, we focus only on the subsets $A_i$ of width at most, say, $2M$.
   Clearly, as $A\subset[n]$, at least $N/2$ of these sets have, indeed, the width at most $2n/N=2M$. Let
   $$I=\{0\le i\le N-1:\; A_i \text{ has width at most } 2M\}.$$
   So, $|I|\ge N/2$. Next, define $2r-1$ disjoint sub-permutations $\pi^{(j)}$ of $\pi$, $j=1,\dots, 2r-1$,  of length $|I|$, by including in $\pi^{(j)}$ all elements of $\pi$ of the form $a_{i(2r-1)+j}$, where $i\in I$ (in the order they appear in $\pi$). To place them on a common ground, we replace each element $a_{i(2r-1)+j}$ by $i$, obtaining a set of new permutations $\bar\pi^{(j)}$, $j=1,\dots, 2r-1$, of the same set $I$.

   We  apply Lemma \ref{3to2r-1} to  permutations $\bar\pi^{(j)}$, $j=1,\dots, 2r-1$, and, as a result, find $r$ of them with a common sub-permutation of length at least
   $$|I|^{1/R}\ge(N/2)^{1/R}=2^{-1/R} n^{1/(2R-1)}\ge2^{-1/3}\frac{N}{M},$$
where the latter inequality follows from $R=\binom{2r-1}{r}\ge \binom{3}{2}=3$.
   Let $\bar\pi^{(j_1)},\dots, \bar\pi^{(j_r)}$, where $1\le j_1<\cdots<j_r\le 2r-1$, be these $r$ permutations and
   let $I_0\subset I\subset\{0,\dots, N-1\}$ be the set of the elements of a common sub-permutation in $\bar\pi^{(j_1)},\dots, \bar\pi^{(j_r)}$ of length $k=|I_0|=\lceil2^{-1/3}\tfrac{N}{M}\rceil$. Upon returning to permutations  $\pi^{(j_1)},\dots, \pi^{(j_r)}$, it can be seen, guided by the above example, that $I_0$ generates in their union, and thus in~$\pi$, a $2M$-narrow $r$-tuplet of length $k$.

   It remains to incorporate the above procedure into an iteration loop and make sure that the obtained $r$-tuplets can be concatenated into a long one.
   Let $\pi$ be a permutation of $[n]$. In the first step of the procedure, take $A^1$ to be the set of the first $(2r-1)N$ elements of $\pi$ and let $\pi^1$ be the sub-permutation of $\pi$ consisting of the first $|A^1|$ elements. By the above described argument we find a $2M$-narrow $r$-tuplet $(x_1^{(j)},\dots,x_k^{(j)})$, $j=1,\dots,r$  of length $k$ in $\pi^1$. We then remove from $[n]$ the set
   $$A^1\cup\bigcup_{i=1}^k[\min_{1\le j\le r} x_i^{(j)},\min_{1\le j\le r} x_i^{(j)}+2M-1].$$
   Note that this set has size at most $(2r-1)N+k(2M)\le (2r+1)N$.

   Now, we consider the truncation of $\pi$ to the remaining elements and the set $A^2$ of its leftmost $(2r-1)N$ elements. Repeating the above procedure mutatis mutandis, we obtain another $2M$-narrow $r$-tuplet $(y_1^{(1)},\dots,y_k^{(r)})$, which, owing to the truncation, satisfies the conditions for proper concatenation with $(x_1^{(j)},\dots,x_k^{(j)})$, $j=1,\dots,r$, spelled out in Subsection~\ref{subsec:prep}. This procedure can be continued as long as there are at least $(2r-1)N$ elements left. Thus, it can be repeated at least $\tfrac n{(2r+1)N}=\tfrac M{2r+1}$ times, yielding at the conclusion an $r$-tuplet in $\pi$ of length a least
   $$\frac M{2r+1}\times k \ge \frac M{2r+1}\times 2^{-1/3}\frac{N}{M} \ge \frac N{3r+1}.$$
   (For $r\ge3$, the last bound can be improved to  $\tfrac N{3r}$.)

\section{Proof of  Theorem \ref{rg1}}

 \subsection{Upper bound}
	For the upper bound we use the first moment method. Let $\Pi=\Pi_n$ be a random permutation chosen uniformly from the set of all $n!$ permutations of $[n]$. Let $k$ be a fixed positive integer and let $X_k$ be a random variable counting all $r$-tuplets of length $k$ in $\Pi$. Furthermore, for a family of $r$ pairwise disjoint subsets $A_1,\dots,A_r$ of $[n]$, each of length $k$, let $X_{A_1,\dots,A_r}$ be an indicator random variable equal to 1 if there is an $r$-tuplet in $\Pi$ on positions determined by the subsets $A_1,\dots,A_r$. So, $X_k=\sum_{A_1,\dots,A_r} X_{A_1,\dots,A_r}$ and by the linearity of expectation $$\EE X_k=\sum_{A_1,\dots,A_r}\EE X_{A_1,\dots,A_r}=\sum_{A_1,\dots,A_r}\PP(X_{A_1,\dots,A_r}=1).$$ Since
\begin{equation*}\label{1|k!}
\PP(X_{A_1,\dots,A_r}=1)=\frac{\binom nk\binom{n-k}k\cdots\binom{n-(r-2)k}k\cdot (n-(r-1)k)!\cdot 1}{n!}=\frac{1}{k!^{r-1}}
\end{equation*}
 and the number of unordered $r$-tuples $\{A_1,\dots, A_r\}$ in $[n]$ is
 $$\binom{n}{rk} \frac{(rk)!}{k!^r r!} =\frac{n!}{k!^r r!(n-rk)!},$$
 it follows that
  $$\EE X_k=\frac{n!}{k!^r r!(n-rk)!} \times \frac{1}{k!^{r-1}}=\frac{n(n-1)\dots (n-rk+1)}{r!(k!)^{2r-1}}.$$ Using the inequality $k!>\frac{k^k}{e^k}$ we obtain
  \begin{equation*}\label{<1}
  \EE X_k<\frac{n^{rk}e^{(2r-1)k}}{r!k^{(2r-1)k}}<\left(\frac{n^re^{2r-1}}{k^{2r-1}}\right)^k.
   \end{equation*}
 It follows that for  $k \ge 2en^{r/(2r-1)}$,
$$\PP(X_k\ge1)\le\E X_k<2^{-(2r-1)k}\to0$$ as $n\to\infty$.
    This completes the proof of the upper bound.


\subsection{Lower bound}
The proof of the lower bound presented here is similar to that in \cite{DGR} (see also \cite{BukhR} for its continuous version).
Set
\begin{equation}\label{a}
a=r!^{1/(2r-1)}n^{(r-1)/(2r-1)},
\end{equation}
assume for convenience that $a$ is an integer and divides $n$, and partition $[n]$ into $n/a$ consecutive blocks of equal size, that is, set
$$[n]=A_1\cup\cdots\cup A_{n/a},$$
where $A_1=\{1,\dots,a\}$, $A_2=\{a+1,\dots,2a\}$, etc.

For fixed $1\le i,j\le n/a$,
let $X=X_{ij}$ be the number of elements from the set $A_j$ which $\Pi$ puts on the positions belonging to the set $A_i$. We construct an auxiliary $n/a\times n/a$ bipartite graph $B$ with vertex classes $U=\{1,\dots,n/a\}$ and $V=\{1,\dots, {n/a}\}$, where $ij\in B$ if, and only if, $X_{ij}\ge r$.

 Let $M=\{i_1j_1,\dots, i_mj_m\}$, $i_1<\cdots<i_m$, be a matching in $B$ of size $|M|=m$. For every $ij\in M$, let $s^{(i)}_1,s^{(i)}_2,\dots,s^{(i)}_r$ be some $r$ elements of $A_i$ such that $\Pi(s^{(i)}_1),\dots,\Pi(s^{(i)}_r)\in A_j$. Then, sub-permutations
 $$(\Pi(s^{(i_1)}_1),\dots,\Pi(s^{(i_m)}_1)),\ (\Pi(s^{(i_1)}_2),\dots,\Pi(s^{(i_m)}_2)),\dots,\;
(\Pi(s^{(i_1)}_r),\dots,\Pi(s^{(i_m)}_r))$$ form an $r$-tuplet. Indeed, if, say, $\Pi(s^{(i_1)}_1)<\Pi(s^{(i_2)}_1)$, then $j_1<j_2$, and so, $\Pi(s^{(i_1)}_t)<\Pi(s^{(i_2)}_t)$ for each $t=2,\dots,r$, as $\Pi(s^{(i_1)}_t)\in A_{j_1}$, while $\Pi(s^{(i_2)}_t)\in A_{j_2}$.

Hence, it remains to show that a.a.s.~there is a matching in $B$ of size $\Omega\left(n^{r/(2r-1)}\right)$. To this end, we are going to use the obvious fact, coming from the greedy algorithm, that in every graph $G$ there is a matching of size at least
$|E(G)|/(2\Delta(G))$,
where $\Delta(G)$ is the maximum vertex degree in $G$.

 In fact, we apply this bound  to a suitably chosen subgraph $B'$ of $B$.
 Let $\nu(B)$ be the size of a largest matching in $B$. Note that $\Delta(B)\le\Delta_0:= \lfloor a/r\rfloor$. Further, let $Z_k$ be the number of vertices of degree $k$ in $B$, $k=0,\dots,\Delta_0$. Then, for the subgraph $B'$ of $B$ obtained by deleting all vertices of degree at least $7$, we get
 $$|E(B')|\ge|E(B)|-\sum_{k=7}^{\Delta_0}kZ_k$$ and
\begin{equation}\label{niu}
\nu(B)\ge\nu(B')\ge\frac{|E(B')|}{2\Delta(B')}\ge\frac{|E(B)|-\sum_{k=7}^{\Delta_0}kZ_k}{12}.
\end{equation}
Further, for each $i\in U\cup V$, let
$$Y_i=\sum_{j=1}^{n/a}\mathbb{I}(X_{ij}\ge r)$$ be the degree of vertex $i$ in $B$. Then,
$$|E(B)|=\frac12\sum_{i\in U}\E Y_i$$
and
$$Z_k=\sum_{i\in U\cup V}\mathbb{I}(Y_i=k),$$
and so
$$\E[|E(B)|]=\frac n{a}\E Y_1$$
and
$$\E Z_k=\frac{2n}a\PP(Y_1=k).$$
Taking the expectation on the outer sides of \eqref{niu}, we arrive at
\begin{equation}\label{EnuB}
\E[\nu(B)]\ge\frac n{12a}\left(\E Y_1-2\sum_{k=7}^{\Delta_0}k\PP(Y_1=k)\right)
\end{equation}
and it remains to estimate $\E Y_1$ and $\PP(Y_1=k)$.

We have
$$\E Y_1=\sum_{j=1}^{n/a}\PP(X_{1j}\ge r)\ge\frac na\PP(X_{11}=r)=\frac na\cdot\frac{\binom{n-a}{a-r}\binom ar^2r!(a-r)!(n-a)!}{n!}\sim\frac{a^{2r-1}}{r!n^{r-1}}\ge\frac12,$$
for sufficiently large $n$, as $a^2=o(n)$ and, by \eqref{a}, $\tfrac{a^{2r-1}}{n^{r-1}}=r!$.


Now we estimate $\PP(Y_1=k)$ for $k\in \{7,\dots,\Delta_0\}$. Observe that, by using inequalities $\binom{m}{p} \le \frac{m^p}{p!}$ and $\tfrac n{n-rk}\le\exp\{rk/(n-rk)\}$,
\[
\PP(Y_1 = k) \le \binom{\frac{n}{a}}{k} \binom{a}{r}^k \binom{a}{rk} \frac{(rk)! (n-rk)!}{n!}\le \frac{a^{2rk-k}e^{(rk)^2/(n-rk)}}{r!^k n^{rk-k} k!}.
\]
 Consequently, since $(rk)^2\le a^2=o(n)$, $k-1\ge 6$, $(k-1)!>\left(\tfrac{k-1}3\right)^{k-1}$, and by \eqref{a},
\begin{align*}
k\PP(Y_1 = k)
\le \frac{a^{2rk-k}(1+o(1))}{r!^k n^{rk-k} (k-1)!} &\le \frac{2a^{2rk-k}}{r!^k n^{rk-k} ((k-1)/3)^k}
\\ &\le\frac{2a^{2r-1}}{r!n^{r-1}}\left( \frac {3a^{2r-1}}{r!(k-1) n^{r-1} } \right)^{k-1}\le4\cdot 2^{-k}.
\end{align*}
 Thus,
\[
\sum_{k=7}^{\Delta_0}k\PP(Y_1=k)  \le 4\sum_{k=7}^{\infty}2^{-k} = \frac1{16}.
\]

Returning to~\eqref{EnuB}, we conclude that
\[
\E[\nu(B)]\ge\frac{n}{12a}\left(\frac{1}{2} - 2\cdot \frac{1}{16} \right)=\frac1{32}\cdot\frac na=\Theta\left(n^{r/(2r-1)}\right).
\]
Since $t^{(r)}(\Pi_n)\ge\nu(B)$, to complete the proof of Theorem~\ref{rg1} it remains to show that the random variable $h(\Pi)=\nu(B)$ is highly concentrated about its mean.
To this end, we are going to use the Azuma-Hoeffding inequality for random permutations (see, e.g., Lemma 11 in~\cite{FP} or  Section 3.2 in~\cite{McDiarmid98}):
\begin{theorem}\label{azuma}
 Let $h(\pi)$ be a function of $n$-permutations such that, for some constant $c>0$, if permutation $\pi_2$ is obtained from permutation $\pi_1$ by swapping two elements, then $|h(\pi_1)-h(\pi_2)|\le c$.
Then, for every $\eta>0$,
\[
\PP(|h(\Pi_n)-\E[h(\Pi_n)]|\ge \eta)\le 2\exp(-\eta^2/(2c^2n)).
\]
\end{theorem}

To apply it to $h(\Pi)=\nu(B)$,  first note that if $\pi_2$ is obtained from a permutation $\pi_1$ by swapping some two of its  elements, then at most two edges of $B$ can be destroyed, and at most two edges can be created, so $|h(\pi_1)-h(\pi_2)|\le2$.
Consequently, Theorem~\ref{azuma} applied with $c=2$ and $\eta=n^{2r/(4r-1)}$ implies that
\[
\PP\left(|h(\Pi_n)-\E[h(\Pi_n)]|\ge n^{2r/(4r-1)}\right)=o(1).
\]
As, crucially, $n^{2r/(4r-1)}=o(n^{r/(2r-1)})$, this completes the proof of the lower bound in Theorem~\ref{rg1}.

\section{Concluding Remarks}

We conclude the paper with two conjectures and two comments related to the  contents of this paper. The first one is a natural generalization of Gawron's conjecture \cite{Gawron} concerning the asymptotic order of the function $t^{(2)}(n)$.

\begin{conj}\label{Con}
	For every $r\geqslant 2$, we have $t^{(r)}(n) = \Theta(n^{r/(2r-1)})$.
\end{conj}

This means that we believe that every permutation contains an $r$-tuplet of length  $\Omega(n^{r/(2r-1)})$. While this statement is wide open for every $r\geqslant 2$, notice that, as $r$ tends to $\infty$, the difference of the exponents in the upper and lower bounds in Theorem \ref{thm:_ulbounds} tends to 0 (since both converge to $1/2$).

Throughout the paper we have kept $r$ fixed, that is, independent of $n$. However, in the next conjecture we allow $r=r(n)$.
We call an $r$-tuplet of length $r$ an \emph{$r$-square}.
The question is how large a square is contained in every  permutation.
Let $s(\pi)$ denote the largest integer $r$ such that $\pi$ contains an $r$-square and let $s(n)$ be the minimum of $s(\pi)$ over all $n$-permutations.
\begin{conj}
	We have $n-(s(n))^2=o(n)$.
\end{conj}
In other words, we believe that every permutation is almost entirely filled by a square.
Notice that a repeated application of the Erd\H {o}s-Szekeres Theorem (\cite{ErdosSzekeres}) implies that every $n$-permutation contains an $\Omega(\sqrt n)$-square. Clearly, the components of this square are monotone. It is, therefore, natural to expect even larger squares with components of an arbitrary pattern.

\smallskip

Next, let us mention $r$-tuplets with forbidden patterns. For a given permutation $\tau$ we say that a permutation $\pi$ is \emph{$\tau$-free} if no subsequence of $\pi$ is similar to $\tau$. For $r\ge2$, let  $t^{(r)}(n,\tau)$ denote the longest length of a $\tau$-free $r$-tuplet guaranteed in every $n$-permutation. In \cite{DGR} we showed that $t^{(2)}(n,\tau)=\Theta(\sqrt n)$ for any non-monotone pattern $\tau$.  It turns out that, somewhat surprisingly,  the same threshold remains valid for any $r\ge2$ (with  only the constant hidden in $\Theta$ possibly depending on $r$). We omit the details.

\medskip

Finally, it should be remarked that there has been an extensive research on twins in other discrete structures, like graphs and words over finite alphabets (see \cite{DGR} for more details and references). Here we would like to point to one of the most interesting results in this area. In 2012, Axenovich, Person, and Puzynina proved in \cite{APP} that every binary sequence of length $n$ contains $r$ disjoint identical subsequences, each of length $n/r-o(n)$. Equally interesting is their method of proof: they proved and then utilized a special version of the regularity lemma for words, an analog of the celebrated Szemer\'edi's Regularity Lemma for graphs.

In the context of the results presented in this paper we wonder if a similar tool would help to improve bounds on $t^{(r)}(n)$ and ultimately prove Conjecture \ref{Con}. We are only aware of a regularity lemma for permutations developed by Cooper in \cite{Coo} (see also \cite{HKS}), but it seems not suitable for tackling the problem of twins.

\subsection*{Acknowledgements} We are very grateful to both referees, as well as to the editors of the journal, for their careful
reading of the manuscript and several  comments which have led to an improvement in clarity and readability of our paper.

\end{document}